\PassOptionsToPackage{unicode}{hyperref}
\PassOptionsToPackage{hyphens}{url}
\PassOptionsToPackage{dvipsnames,svgnames,x11names}{xcolor}
\documentclass[
  11pt,
  a4paper,
]{article}
\usepackage{amsmath,amssymb}
\usepackage{iftex}
\ifPDFTeX
  \usepackage[T1]{fontenc}
  \usepackage[utf8]{inputenc}
  \usepackage{textcomp} 
\else 
  \usepackage{unicode-math} 
  \defaultfontfeatures{Scale=MatchLowercase}
  \defaultfontfeatures[\rmfamily]{Ligatures=TeX,Scale=1}
\fi
\usepackage{lmodern}
\ifPDFTeX\else
\fi
\IfFileExists{upquote.sty}{\usepackage{upquote}}{}
\IfFileExists{microtype.sty}{
  \usepackage[]{microtype}
  \UseMicrotypeSet[protrusion]{basicmath} 
}{}
\makeatletter
\@ifundefined{KOMAClassName}{
  \IfFileExists{parskip.sty}{%
    \usepackage{parskip}
  }{
    \setlength{\parindent}{0pt}
    \setlength{\parskip}{6pt plus 2pt minus 1pt}}
}{
  \KOMAoptions{parskip=half}}
\makeatother
\usepackage{xcolor}
\usepackage[margin=1in]{geometry}
\usepackage{longtable,booktabs,array}
\usepackage{calc} 
\usepackage{etoolbox}
\makeatletter
\patchcmd\longtable{\par}{\if@noskipsec\mbox{}\fi\par}{}{}
\makeatother
\IfFileExists{footnotehyper.sty}{\usepackage{footnotehyper}}{\usepackage{footnote}}
\makesavenoteenv{longtable}
\providecommand{\tightlist}{%
  \setlength{\itemsep}{0pt}\setlength{\parskip}{0pt}}
\ifLuaTeX
  \usepackage{selnolig}  
\fi
\IfFileExists{bookmark.sty}{\usepackage{bookmark}}{\usepackage{hyperref}}
\IfFileExists{xurl.sty}{\usepackage{xurl}}{} 
\hypersetup{
  pdftitle={Finite and Dynamic Stability Horizons for Nearest-Neighbor Future Structures},
  pdfauthor={Hideki Sumiya},
  colorlinks=true,
  linkcolor={blue},
  filecolor={Maroon},
  citecolor={Blue},
  urlcolor={blue},
  pdfcreator={LaTeX via pandoc}}

\title{Finite and Dynamic Stability Horizons for Nearest-Neighbor Future Structures}
\author{Hideki Sumiya\\
\small Department of Orthopaedic Surgery, Nozaki Tokushukai Clinic,\\
\small Daito, Osaka, Japan\\
\small \texttt{puchanorange@icloud.com}}
\date{}

\begin{document}
\maketitle

\begin{abstract}

Nearest-neighbor graphs are discrete objects whose membership may change
under small perturbations of the underlying coordinates. We establish an
explicit stability guarantee for finite labeled configurations in an
arbitrary metric space. If the gap between the \(k\)th and \((k+1)\)st
distances is positive, simultaneous per-label perturbations smaller than
one quarter of that gap preserve the directed \(k\)-nearest-neighbor
membership. The factor four is sharp under the stated uniform
displacement assumptions, and every fixed deterministic construction
based solely on the labeled neighbor family is consequently invariant.
Under an interval-valid Lipschitz bound for labeled information-space
trajectories, the same result yields a certified lower bound on the
first possible rewiring time, with a refinement for label-specific
motion bounds. We apply the finite theorem to frozen standardized
Taylor--Green future-information coordinates \([d_B,\log A_B]\) for
585 particles at three observed time strata. Outward-rounded interval
arithmetic certified a sufficient perturbation radius of approximately
\(1.023\times10^{-6}\), and a separately implemented checker within the
same research workflow verified the rank and distance-margin calculations.
An outcome-blind audit then examined
whether the construction supported a numerical continuous-time horizon.
Because the complete information map involved discrete clustering and
boundary reconstruction and lacked an analytic derivative bound,
validated dense-time bound, or certified modulus of continuity, the
application was correctly classified as \texttt{DISCRETE\_ONLY} with
\texttt{STOP\_NO\_INTERVAL\_VALID\_BOUND}. Thus finite local structural
invariance is proved and numerically certified for the frozen
configuration, while temporal specialization and predictive
generalization remain separate questions requiring additional evidence.
\end{abstract}


\hypertarget{introduction}{%
\section{1. Introduction}\label{introduction}}

Nearest-neighbor rules and neighborhood graphs are foundational tools in
nonparametric classification and in geometric approaches to
dimensionality reduction and data representation {[}1--3{]}. In these
settings, a basic interpretive question is whether the observed
adjacency is a stable feature of the configuration or an artifact of
small coordinate perturbations. Empirical persistence across a few
samples can be suggestive, but it does not by itself provide a
perturbation certificate.

The relevant obstruction is rank degeneracy. For a directed
\(k\)-nearest-neighbor rule, membership can change only when a selected
and an unselected candidate exchange their distance order. The gap
between the \(k\)th and \((k+1)\)st distances therefore measures separation
from the nearest membership-changing boundary. This observation is
elementary, but making it operational requires care about simultaneous
movement of the focal point and both competing labels, strict
inequalities, tie conventions, and deterministic constructions
downstream of the graph.

Robustness and certification have also been studied for nearest-neighbor
classifiers under adversarial examples, uncertain data, and abstract
perturbation models {[}4--6{]}. Those studies concern classification
outcomes, test-point perturbations, or alternative admissible datasets.
The problem considered here is different: every label in a finite metric
configuration may move simultaneously, and the object to be certified is
the complete labeled directed neighbor family together with any fixed
deterministic graph-based cascade. We therefore do not claim that
robustness of k-NN methods is generally new; we isolate and solve this
specific whole-configuration invariance problem.

Temporal stability is a separate question. A positive rank margin at one
time gives a spatial perturbation budget. Converting that budget into a
time horizon additionally requires an interval-valid bound on the motion
of the same information-space coordinates. Endpoint differences or
observed-time secant slopes do not supply such a bound. This distinction
becomes especially important when an information coordinate is rebuilt
through clustering, boundary extraction, or other discrete operations,
because smooth motion of the underlying physical state need not imply a
certified smooth trajectory for the derived coordinate.

This study develops a finite and dynamic stability framework for labeled
directed \(k\)-nearest-neighbor structures. First, we prove in an
arbitrary metric space that a positive selected--unselected boundary
margin yields an explicit uniform perturbation radius preserving every
directed neighbor set. The factor four in the bound accounts for
simultaneous motion of the focal, selected, and unselected labels and is
sharp under the stated assumptions. Second, we show that every fixed
deterministic construction based only on the labeled neighbor family is
exactly invariant inside that radius. Third, under an interval-valid
Lipschitz hypothesis, we convert the static margin into a certified
lower bound on the first possible rewiring time, including a
label-specific refinement.

We then specialize the finite theorem to a frozen Taylor--Green
future-information construction. At three observed time strata, 585
particles were embedded in standardized coordinates \([d_B,\log A_B]\), and directed 5-nearest-neighbor graphs were formed within each
stratum. Exact outward-rounded interval calculations certified a
positive global boundary margin and hence a sufficient perturbation
radius of approximately \(1.023\times10^{-6}\) in the frozen
standardized Euclidean coordinates. A separately implemented checker within
the same research workflow verified the rank identities, distance enclosures,
and margin identities. Because the
downstream graph statistics and labels were deterministic functions of
the frozen neighbor families, they inherit exact invariance inside the
certified radius.

Finally, we audit whether the same construction supports a numerical
continuous-time stability horizon. The answer is presently negative for
a precise reason: although the underlying trajectories and deformation
gradients were evaluated densely, the complete information map was
reconstructed only at discrete times and included clustering and
boundary-selection operations. No analytic derivative bound, validated
dense-time bound, or certified modulus of continuity was available for
\([d_B,\log A_B]\). The correct specialization is therefore static at
the observed strata; the dynamic theorem remains conditional and valid,
but its Taylor--Green time-horizon hypothesis was not established.

The contribution is thus a separation of three claims that are often
conflated: observed persistence, certified local perturbation
invariance, and certified temporal invariance. The first is empirical,
the second follows here from a finite rank-margin certificate, and the
third requires additional interval control. None of these claims alone
establishes predictive validity, causality, stability of the underlying
Navier--Stokes trajectories, or universality across flows. The
prospective F32B-T analysis, in particular, did not support the complete
preregistered phase-switch pattern; that result is reported as
scientifically distinct from the mathematical stability certificate.

\hypertarget{mathematical-setting}{%
\section{2. Mathematical setting}\label{mathematical-setting}}

Let \((X,d)\) be a metric space, let \(n\ge3\), and let

\[Z=(z_1,\ldots,z_n)\in X^n\]

be a finite labeled configuration. Fix \(k\) with \(1\le k\le n-2\). For
each focal label \(i\), order the distances to the remaining labels
using the lexicographic rule \((d(z_i,z_j),j)\):

\[d_{i,(1)}(Z)\le\cdots\le d_{i,(n-1)}(Z).\]

The directed \(k\)-nearest-neighbor set is denoted by \(N_i^{(k)}(Z)\).
Its boundary margin is

\[\gamma_i(Z)=d_{i,(k+1)}(Z)-d_{i,(k)}(Z),\]

and the global margin is

\[\gamma_{\min}(Z)=\min_{1\le i\le n}\gamma_i(Z).\]

A positive boundary margin excludes a tie across the
selected--unselected boundary. Ties within either side do not affect
membership. For two labeled configurations define

\[D_\infty(Z,\widetilde Z)=\max_i d(z_i,\widetilde z_i).\]

The case \(k=n-1\) is excluded from the margin definition because all
other labels are selected; its neighbor family is trivially invariant
under every label-preserving perturbation.

\hypertarget{finite-stability-of-future-structures}{%
\section{3. Finite stability of future
structures}\label{finite-stability-of-future-structures}}

\hypertarget{theorem-1-finite-metric-k-nn-stability}{%
\subsection{Theorem 1 --- finite metric k-NN
stability}\label{theorem-1-finite-metric-k-nn-stability}}

Suppose \(\gamma_i(Z)>0\) for a fixed focal label \(i\). If

\[D_\infty(Z,\widetilde Z)<\frac{\gamma_i(Z)}4,\]

then

\[N_i^{(k)}(\widetilde Z)=N_i^{(k)}(Z).\]

Consequently, if \(\gamma_{\min}(Z)>0\) and

\[D_\infty(Z,\widetilde Z)<\frac{\gamma_{\min}(Z)}4,\]

then

\[N_i^{(k)}(\widetilde Z)=N_i^{(k)}(Z)\qquad\text{for every }i.\]

\hypertarget{proof}{%
\subsection{Proof}\label{proof}}

Put \(r=D_\infty(Z,\widetilde Z)\). The triangle inequality implies

\[|d(\widetilde z_i,\widetilde z_j)-d(z_i,z_j)|\le2r.\]

For \(j\in N_i^{(k)}(Z)\) and \(\ell\notin N_i^{(k)}(Z)\),

\[d(z_i,z_\ell)-d(z_i,z_j)\ge\gamma_i(Z).\]

Hence

\[d(\widetilde z_i,\widetilde z_\ell)-d(\widetilde z_i,\widetilde z_j)\ge\gamma_i(Z)-4r>0.\]

Every initially selected label therefore remains closer than every
initially unselected label for the fixed focal label \(i\). Applying the
same argument to every \(i\) under the global margin condition proves
full graph invariance. \(\square\)

\hypertarget{corollary-1-deterministic-cascade-invariance}{%
\subsection{Corollary 1 --- deterministic cascade
invariance}\label{corollary-1-deterministic-cascade-invariance}}

If \(\mathcal A(Z)=F(N_1^{(k)}(Z),\ldots,N_n^{(k)}(Z))\) for a fixed
deterministic rule \(F\), then

\[
\gamma_{\min}(Z)>0,
\qquad
D_\infty(Z,\widetilde Z)<\frac{\gamma_{\min}(Z)}4
\]

imply \(\mathcal A(\widetilde Z)=\mathcal A(Z)\).

\hypertarget{corollary-2-necessary-displacement-for-rewiring}{%
\subsection{Corollary 2 --- necessary displacement for
rewiring}\label{corollary-2-necessary-displacement-for-rewiring}}

If \(N_i^{(k)}(\widetilde Z)\ne N_i^{(k)}(Z)\), then

\[D_\infty(Z,\widetilde Z)\ge\frac{\gamma_i(Z)}4.\]

\hypertarget{proposition-1-sharpness-of-the-universal-factor-four}{%
\subsection{Proposition 1 --- sharpness of the universal factor
four}\label{proposition-1-sharpness-of-the-universal-factor-four}}

No constant \(C<4\) can replace \(4\) in Theorem 1 under only a uniform
per-label displacement bound.

\hypertarget{proof-1}{%
\subsection{Proof}\label{proof-1}}

Fix \(C<4\). On \(\mathbb R\), take \(n=3\), \(k=1\), a focal point at
\(0\), a selected point at \(-a\), and an unselected point at
\(a+\gamma\), where \(a,\gamma>0\). The initial boundary margin is
\(\gamma\). Under the perturbation

\[
0\mapsto\varepsilon,\qquad
-a\mapsto-a-\varepsilon,\qquad
a+\gamma\mapsto a+\gamma-\varepsilon,
\]

the selected and unselected distances from the focal point become
\(a+2\varepsilon\) and \(a+\gamma-2\varepsilon\), respectively. Thus
their gap is \(\gamma-4\varepsilon\). Choose

\[
\frac{\gamma}{4}<\varepsilon<\frac{\gamma}{C},
\]

which is possible because \(C<4\). Then
\(D_\infty=\varepsilon<\gamma/C\), yet the nearest-neighbor membership
reverses strictly. Hence no universal denominator smaller than four is
valid. At \(\varepsilon=\gamma/4\), a boundary tie occurs, showing why
the strict inequality in Theorem 1 is necessary. \(\square\)

\hypertarget{dynamic-stability-horizons}{%
\section{4. Dynamic stability
horizons}\label{dynamic-stability-horizons}}

Let \(Z(t)=(z_1(t),\ldots,z_n(t))\) be labeled trajectories on
\([t_0,t_0+T]\), with fixed metric, labels, value of \(k\), and
tie-breaking rule.

\hypertarget{theorem-2-uniform-dynamic-stability-horizon}{%
\subsection{Theorem 2 --- uniform dynamic stability
horizon}\label{theorem-2-uniform-dynamic-stability-horizon}}

Assume every trajectory is \(L\)-Lipschitz:

\[d(z_a(t),z_a(u))\le L|t-u|.\]

If \(L>0\), define

\[\tau_i^{\mathrm{cert}}(t_0)=\min\left\{T,\frac{\gamma_i(Z(t_0))}{4L}\right\}.\]

Then

\[N_i^{(k)}(Z(t_0+s))=N_i^{(k)}(Z(t_0))\]

for every \(0\le s<\tau_i^{\mathrm{cert}}(t_0)\). If
\(T<\gamma_i(Z(t_0))/(4L)\), invariance also holds at \(s=T\). If \(L=0\), it
holds throughout \([0,T]\).

\hypertarget{proof-2}{%
\subsection{Proof}\label{proof-2}}

Lipschitz continuity gives \(D_\infty(Z(t_0),Z(t_0+s))\le Ls\). The
claim follows from the focal-label clause of Theorem 1. \(\square\)

\hypertarget{corollary-3-first-rewiring-lower-bound}{%
\subsection{Corollary 3 --- first-rewiring lower
bound}\label{corollary-3-first-rewiring-lower-bound}}

With \(\inf\varnothing=+\infty\), define

\[\tau_i^{\mathrm{rewire}}=\inf\{s\in(0,T]:N_i^{(k)}(Z(t_0+s))\ne N_i^{(k)}(Z(t_0))\}.\]

Then, whenever \(L>0\),

\[\tau_i^{\mathrm{rewire}}\ge\frac{\gamma_i(Z(t_0))}{4L}.\]

Here the infimum is taken only over \((0,T]\); if the displayed lower
bound exceeds \(T\), Theorem 2 implies that the set is empty and the
adopted value is \(+\infty\).

\hypertarget{theorem-3-label-specific-refinement}{%
\subsection{Theorem 3 --- label-specific
refinement}\label{theorem-3-label-specific-refinement}}

Suppose \(d(z_a(t_0+s),z_a(t_0))\le L_as\) for every label and
\(s\in[0,T]\). For selected \(j\) and unselected \(\ell\), put

\[g_{ij\ell}=d(z_i(t_0),z_\ell(t_0))-d(z_i(t_0),z_j(t_0)).\]

If \(\gamma_i(Z(t_0))>0\), define

\[
\tau_{i,\mathrm{loc}}^{\mathrm{cert}}
=
\min\left\{
T,
\min_{j\in N_i^{(k)},\,\ell\notin N_i^{(k)}}
\frac{g_{ij\ell}}{2L_i+L_j+L_\ell}
\right\}.
\]

Here a zero denominator is interpreted as \(+\infty\). The neighbor set
is invariant for every \(0\le s<\tau_{i,\mathrm{loc}}^{\mathrm{cert}}\).
If \(T\) is strictly smaller than the inner minimum, invariance also
holds at \(s=T\).

\hypertarget{proof-3}{%
\subsection{Proof}\label{proof-3}}

Fix \(j\in N_i^{(k)}(Z(t_0))\), \(\ell\notin N_i^{(k)}(Z(t_0))\),
and \(s\in[0,T]\). The triangle inequality and the assumed
displacement bounds give

\[
d(z_i(t_0+s),z_\ell(t_0+s))
\ge d(z_i(t_0),z_\ell(t_0))-(L_i+L_\ell)s
\]

and

\[
d(z_i(t_0+s),z_j(t_0+s))
\le d(z_i(t_0),z_j(t_0))+(L_i+L_j)s.
\]

Subtracting yields

\[
d(z_i(t_0+s),z_\ell(t_0+s))
-d(z_i(t_0+s),z_j(t_0+s))
\ge g_{ij\ell}-(2L_i+L_j+L_\ell)s.
\]

For \(s<\tau_{i,\mathrm{loc}}^{\mathrm{cert}}\), the right-hand side
is strictly positive for every selected--unselected pair with nonzero
denominator; when the denominator is zero, the same bound leaves the
initial positive gap unchanged. Hence every initially selected label
remains closer than every initially unselected label. The endpoint
claim follows identically when \(T\) is strictly smaller than the inner
minimum. \(\square\)

\hypertarget{remark-applicability}{%
\subsection{Remark --- applicability}\label{remark-applicability}}

The dynamic theorem requires an interval-valid bound for the chosen
information map. Endpoint secant slopes do not provide such a bound. In
the current Taylor--Green construction, \(d_B\) is rebuilt at discrete
times through clustering and boundary edges; consequently F34 permits
only static certificates at the observed time strata and does not
authorize a numerical continuous-time horizon, as detailed in Section 6.

\hypertarget{certified-finite-stability-in-a-taylorgreen-future-information-graph}{%
\section{5. Certified finite stability in a Taylor--Green
future-information
graph}\label{certified-finite-stability-in-a-taylorgreen-future-information-graph}}

\hypertarget{frozen-information-space-configuration}{%
\subsection{5.1 Frozen information-space
configuration}\label{frozen-information-space-configuration}}

We applied Theorem 1 to the frozen F29 Taylor--Green future-information
construction. The analysis used three separate time strata,

\[t\in\{8.9,9.0,9.1\},\]

each containing the same 585 labeled particles. The 1755 particle--time
records therefore represent three within-time configurations; they do
not form a single nearest-neighbor graph in which points from different
times are mutually compared.

For every stratum, particle \(i\) was represented in the frozen,
globally standardized two-dimensional information space

\[z_i(t)=\bigl[d_{B,i}(t),\log A_{B,i}(t)\bigr]\in\mathbb R^2.\]

Euclidean distance and directed \(k=5\) nearest-neighbor membership were
used. The undirected information graph at each time was the symmetrized
union of its directed 5-NN edges. All normalization parameters, particle
labels, time keys, distance conventions, and downstream rules were held
fixed.

\hypertarget{certified-rank-margin}{%
\subsection{5.2 Certified rank margin}\label{certified-rank-margin}}

For each particle and time stratum, the directed boundary margin was

\[\gamma_i(t)=d_{i,(6)}(t)-d_{i,(5)}(t),\]

and the global stratified margin was

\[\gamma_{\min,T}=\min_{t\in\{8.9,9.0,9.1\}}\min_i\gamma_i(t).\]

Outward-rounded interval calculations certified the strict lower bound

\[
\gamma_{\min,T}>4.0919760\times10^{-6}.
\]

Consequently, a conservative decimal certificate for Theorem 1 is

\[
r_{\mathrm{cert}}=1.0229940\times10^{-6}.
\]

The complete outward-rounded interval endpoints are retained in the
machine-readable certificate rather than reported as physically
meaningful decimal precision.

For every simultaneous perturbation of the frozen standardized
particle--time coordinates satisfying

\[
\max_{t,i}\|\widetilde z_i(t)-z_i(t)\|_2
<r_{\mathrm{cert}},
\]

all directed 5-NN memberships are unchanged at all three time strata.
The strict inequality is essential. The minimizing observation was
particle 901 at \(t=9.1\).

\hypertarget{exact-downstream-consequences}{%
\subsection{5.3 Exact downstream
consequences}\label{exact-downstream-consequences}}

Because the labeled directed neighbor families are invariant inside the
certified radius, every downstream object computed deterministically
from those families and the frozen conventions is also invariant. This
statement requires that no additional coordinate-dependent, refitted,
randomized, or externally perturbed input enter after graph
construction. Under that dependency condition, the invariant objects in
the F30 construction were:

\begin{enumerate}
\def\labelenumi{\arabic{enumi}.}
\tightlist
\item
  directed 5-NN memberships;
\item
  symmetrized undirected information graphs \(G_t\);
\item
  undirected local clustering coefficients \(C_i(t)\);
\item
  within-time NumPy-linear 90th-percentile thresholds \(Q90_t\);
\item
  strict-core labels \(\mathbf 1\{C_i(t)>Q90_t\}\);
\item
  core-neighbor exposures \(E_i(t)\);
\item
  observational next-time strict-core-entry indicators \(Y_i\).
\end{enumerate}

No additional Lipschitz constant is required for this deterministic
cascade. The conclusion follows because the labeled discrete inputs to
every downstream rule are identical. This invariance does not establish
that the downstream quantities predict a physical transition.

\hypertarget{second-implementation-numerical-verification}{%
\subsection{5.4 Second-implementation numerical
verification}\label{second-implementation-numerical-verification}}

The separately implemented F30 checker, developed within the same research
workflow rather than by an external reviewer, reported PASS after verifying:

\begin{itemize}
\tightlist
\item
  1755 rank identities;
\item
  3510 squared distances;
\item
  3510 outward square-root enclosures;
\item
  1755 margin identities.
\end{itemize}

The principal frozen records were:

\begin{longtable}[]{@{}
  >{\raggedright\arraybackslash}p{(\columnwidth - 2\tabcolsep) * \real{0.5000}}
  >{\raggedright\arraybackslash}p{(\columnwidth - 2\tabcolsep) * \real{0.5000}}@{}}
\toprule\noalign{}
\begin{minipage}[b]{\linewidth}\raggedright
Artifact
\end{minipage} & \begin{minipage}[b]{\linewidth}\raggedright
SHA-256
\end{minipage} \\
\midrule\noalign{}
\endhead
\bottomrule\noalign{}
\endlastfoot
\path{F30_phase3_finite_stability_theorem_integration.json} &
\url{37dcdcf9c17d72f6c2ba5009bcb9e9fe035f6fa17bfb9ed697f5de7babb11206} \\
\path{F30_phase4_independent_checker_report.json} &
\url{3f1d67f413b3b614411cf4ca7c632caf03bbe7da6471e69792f4ddffa1aa627a} \\
\end{longtable}

The Phase 2 runtime correction affected only canonicalization of time
keys within \(10^{-12}\). The certificate value had not been opened
before correction; the corrected runner was preregistered, the original
runner was retained, and the mathematical and scientific methods were
unchanged.

\hypertarget{interpretation-of-the-certificate}{%
\subsection{5.5 Interpretation of the
certificate}\label{interpretation-of-the-certificate}}

The value \(r_{\mathrm{cert}}\) is a certified sufficient radius based
on a rigorous lower enclosure of \(\gamma_{\min,T}\). It is not asserted
to be the exact maximal configuration-specific robustness radius. The
factor-four theorem uses a worst-case triangle-inequality bound; a
larger safe radius may exist for this particular configuration and would
require a separate constrained optimization proof.

The certificate is expressed in the frozen standardized coordinates. If
normalization is recomputed after perturbation, the induced coordinate
displacement must be bounded and included in the perturbation budget.
The present certificate cannot be transferred automatically to raw or
restandardized coordinates.

\hypertarget{separation-from-temporal-and-predictive-claims}{%
\subsection{5.6 Separation from temporal and predictive
claims}\label{separation-from-temporal-and-predictive-claims}}

The F30 certificate is static at each of the three observed time strata.
It proves local invariance under bounded perturbations of the frozen
configurations. It does not provide an interval-valid temporal speed
bound for

\[\Psi=[d_B,\log A_B],\]

and therefore does not yield a numerical continuous-time stability
horizon. As detailed in Section 6, the outcome-blind F34 audit formally
classified the available Taylor--Green artifacts as
\texttt{DISCRETE\_ONLY}: only static certificates at the observed strata
are authorized.

F30 is also logically distinct from the F32B-T temporal holdout. The
former asks whether the finite information structure is locally
invariant; the latter asks whether a preregistered phase-switch pattern
reproduces prospectively. The valid combined statement is

\[
\boxed{\begin{aligned}
&\text{general conditional finite stability: proved;}\\
&\text{complete prospective phase-switch pattern: not supported.}
\end{aligned}}
\]

The negative F32B-T result does not weaken the deterministic theorem,
and the deterministic theorem does not provide predictive validation.

\hypertarget{exceptional-set-interpretation}{%
\subsection{5.7 Exceptional-set
interpretation}\label{exceptional-set-interpretation}}

Let

\[
\Sigma_k=
\left\{Z:\min_i\bigl[d_{i,(k+1)}(Z)-d_{i,(k)}(Z)\bigr]=0\right\}
\]

denote the k-NN boundary-degeneracy set. Away from \(\Sigma_k\), the
labeled neighbor structure is locally constant. The F30 certificate
identifies a protected neighborhood contained within one combinatorial
cell of the frozen configuration space. Reaching the degeneracy set
permits a change in adjacency but does not by itself prove that a
particular physical or predictive transition occurs.

\hypertarget{outcome-blind-audit-of-continuous-time-applicability}{%
\section{6. Outcome-blind audit of continuous-time
applicability}\label{outcome-blind-audit-of-continuous-time-applicability}}

\hypertarget{question-and-decision-rule}{%
\subsection{6.1 Question and decision
rule}\label{question-and-decision-rule}}

The dynamic stability theorem converts a positive rank margin into a
certified time horizon only when the motion of the chosen
information-space trajectories is bounded throughout an interval. For

\[z_i(t)=\Psi_H(x_i(t))=[d_{B,i}(t),\log A_{B,i}(t)],\]

a continuous-time application therefore requires either an
interval-valid Lipschitz bound

\[\|z_i(t)-z_i(u)\|_2\le L_\Psi|t-u|,\]

or a certified modulus of continuity \(\omega\) satisfying

\[\|z_i(t+s)-z_i(t)\|_2\le\omega(s).\]

F34 was preregistered as an outcome-blind feasibility audit. It was not
permitted to calculate a time horizon unless an analytic derivative
bound, a validated dense-time numerical bound, or a certified modulus
was established. Files containing strict-core outcomes, F31/F32 results,
\(\varepsilon_c\), susceptibility, or intervention-response quantities
were outside the audit.

\hypertarget{available-temporal-computation}{%
\subsection{6.2 Available temporal
computation}\label{available-temporal-computation}}

The underlying velocity field was available as a time-indexed array and
was evaluated using linear interpolation in time and periodic cubic
B-spline interpolation in space. Particle trajectories and deformation
gradients \(J_F\) were integrated by fourth-order Runge--Kutta with

\[\Delta t=0.002,\]

from starting times

\[8.9,\ 9.0,\ 9.1\]

to an end time of \(9.2\). These calculations provided dense evaluations
of the underlying flow map required for \(J_F\). They did not, however,
provide dense evaluations of the complete information map \(\Psi_H\).

The initial F34 preregistration listed \(9.18\) as the fourth known
time, inherited from an earlier analysis convention. Direct
code-provenance inspection subsequently established that the geometry
builders use \(8.9,9.0,9.1,9.2\). The \(J_F\) builders use starting
times \(8.9,9.0,9.1\) and integrate to \(9.2\); the AB builders use the
resulting \(J_F\) values at those three starting times. The
preregistered \(9.18\) entry was therefore a documentation error and was
not an AB observation used in this audit. This correction does not
affect the feasibility decision: neither time list supplies an
interval-valid bound for \(\Psi_H\).

\hypertarget{discrete-reconstruction-of-the-information-coordinates}{%
\subsection{6.3 Discrete reconstruction of the information
coordinates}\label{discrete-reconstruction-of-the-information-coordinates}}

The boundary-distance coordinate \(d_B\) was constructed separately at

\[t\in\{8.9,9.0,9.1,9.2\}\]

from outcome-blind 18-dimensional future-trajectory signatures, cluster
assignments, boundary edges, and boundary geometry. The coordinate
therefore depends on discrete combinatorial operations. Changes in
cluster membership, nearest boundary edge, or boundary connectivity can
change the defining branch of \(d_B\).

The amplification coordinate \(A_B\) was subsequently constructed at the
three starting times by combining the boundary normal with \(J_F\),
followed by the transformation \(\log A_B\). The F28 AB table contained
1755 rows, corresponding to 585 particles at each of three starting
times. The F29 builder used the same 585-by-3 design. Both the F28 and
F29 checks reported that all sampled \(A_B\) values were finite and
positive; for F28, the maximum boundary-normal unit error was

\[2.220446049250313\times10^{-16}.\]

Sampled positivity makes \(\log A_B\) well defined at the recorded rows.
It does not prove positivity or regularity at every intervening time.

\hypertarget{evaluation-of-the-four-certification-routes}{%
\subsection{6.4 Evaluation of the four certification
routes}\label{evaluation-of-the-four-certification-routes}}

\hypertarget{route-a-analytic-derivative-bound}{%
\paragraph{Route A: analytic derivative
bound}\label{route-a-analytic-derivative-bound}}

Route A was not established. No bound was propagated through the entire
chain from the flow field to future signatures, clustering,
boundary-edge selection, \(d_B\), boundary normals, \(J_F\), and
\(\log A_B\). In particular, the discrete changes of cluster and edge
identities were not controlled by a differentiable formula or a
piecewise-smooth switching analysis.

\hypertarget{route-b-validated-dense-time-numerical-bound}{%
\paragraph{Route B: validated dense-time numerical
bound}\label{route-b-validated-dense-time-numerical-bound}}

Route B was not established. Dense integration of trajectories and
\(J_F\) is not equivalent to dense evaluation of \([d_B,\log A_B]\).
Moreover, the available pipeline did not provide a rigorous combined
bound for temporal interpolation, spatial interpolation, Runge--Kutta
truncation, floating-point error, and boundary-reconstruction error in
the transformed coordinates.

\hypertarget{route-c-certified-modulus-of-continuity}{%
\paragraph{Route C: certified modulus of
continuity}\label{route-c-certified-modulus-of-continuity}}

Route C was not established. No modulus bounded the change of \(d_B\)
across possible cluster, edge, or boundary reorganization. Continuity of
the underlying particle trajectories alone does not supply continuity of
this data-defined boundary map.

\hypertarget{route-d-observed-time-static-certificates}{%
\paragraph{Route D: observed-time static
certificates}\label{route-d-observed-time-static-certificates}}

Route D was adopted. The available artifacts support comparisons and
static perturbation certificates at the observed time strata. Endpoint
displacements and secant slopes may be reported descriptively, but they
do not bound the maximum speed between endpoints.

\hypertarget{formal-audit-result}{%
\subsection{6.5 Formal audit result}\label{formal-audit-result}}

The F34 decision, reproduced by a separate implementation within the same
research workflow, was

\[\boxed{\texttt{DISCRETE\_ONLY}}\]

with the continuous-time claim status

\[\boxed{\texttt{STOP\_NO\_INTERVAL\_VALID\_BOUND}}.\]

The audit gates were resolved as follows:

\begin{longtable}[]{@{}
  >{\raggedright\arraybackslash}p{(\columnwidth - 4\tabcolsep) * \real{0.3333}}
  >{\raggedright\arraybackslash}p{(\columnwidth - 4\tabcolsep) * \real{0.3333}}
  >{\raggedright\arraybackslash}p{(\columnwidth - 4\tabcolsep) * \real{0.3333}}@{}}
\toprule\noalign{}
\begin{minipage}[b]{\linewidth}\raggedright
Gate
\end{minipage} & \begin{minipage}[b]{\linewidth}\raggedright
Result
\end{minipage} & \begin{minipage}[b]{\linewidth}\raggedright
Meaning
\end{minipage} \\
\midrule\noalign{}
\endhead
\bottomrule\noalign{}
\endlastfoot
Provenance & PASS for the discrete lineage & AB, geometry, and \(J_F\)
builders were identified and hash-fixed \\
Sampled positivity & PASS & Recorded F28 and F29 \(A_B\) values were
finite and positive \\
Interval regularity & NOT ESTABLISHED & Positivity and smoothness
between sampled times were not certified \\
Temporal validity & FAIL & No interval-valid \(L_\Psi\) or modulus was
established \\
Independence/nonvacuity & NOT REACHED & No time-horizon calculation was
authorized \\
\end{longtable}

Accordingly, the finite-difference quantity

\[
\max_m\frac{\|z_i(t_{m+1})-z_i(t_m)\|_2}{t_{m+1}-t_m}
\]

must not be labeled \(L_\Psi\), and the numerical expression

\[
\tau_{\mathrm{graph}}^{\mathrm{cert}}
=\min\left\{T,\frac{\gamma_{\min}(Z(t_0))}{4L_\Psi}\right\}
\]

must not be evaluated from the current artifacts.

\hypertarget{interpretation}{%
\subsection{6.6 Interpretation}\label{interpretation}}

This STOP does not weaken the static F30 certificate in Section 5 or the
dynamic theorem in Theorem 2. The static theorem applies to fixed
observed configurations and was exactly certified. The dynamic theorem
is a valid conditional statement in any metric space. F34 shows only
that the current Taylor--Green information map has not yet been shown to
satisfy the theorem's temporal hypothesis.

The distinction is substantive:

\[
\begin{aligned}
&\text{dense underlying flow trajectory}\\[-2pt]
&\qquad\not\Longrightarrow
\text{certified dense trajectory of a data-defined boundary coordinate}.
\end{aligned}
\]

Thus the correct conclusion is not that the future-information structure
lacks temporal stability, but that its continuous-time stability horizon
is presently unidentifiable from the frozen construction.

\hypertarget{requirements-for-a-future-continuous-time-specialization}{%
\subsection{6.7 Requirements for a future continuous-time
specialization}\label{requirements-for-a-future-continuous-time-specialization}}

A later study could reopen the question only with a separately
preregistered construction that provides at least one of the following:

\begin{enumerate}
\def\labelenumi{\arabic{enumi}.}
\tightlist
\item
  a fixed smooth information map with an analytic derivative bound;
\item
  a piecewise-smooth boundary map together with certified switching
  times and one-sided bounds;
\item
  dense-time recomputation of the complete boundary geometry with
  validated numerical remainder bounds;
\item
  a certified modulus of continuity that remains valid across boundary
  and cluster changes.
\end{enumerate}

Such a study would constitute new work. It must not retroactively alter
F30 or reinterpret the F34 STOP decision.

\hypertarget{discussion}{%
\section{7. Discussion}\label{discussion}}

\hypertarget{what-the-certificate-establishes}{%
\subsection{7.1 What the certificate
establishes}\label{what-the-certificate-establishes}}

The main mathematical result identifies an explicit protected
neighborhood of any finite labeled configuration with positive directed
\(k\)-NN boundary margins. Within this neighborhood, the graph does not
merely change by a small amount: its labeled directed adjacency is
identical. Any fixed deterministic downstream construction from that
adjacency is therefore identical as well. This is a combinatorial
invariance statement rather than a continuity estimate.

The result can be viewed geometrically. Configuration space is
partitioned into cells on which the labeled neighbor family is constant,
separated by degeneracy loci where selected and unselected distance
ranks meet. A positive boundary margin certifies separation from those
loci. The global radius \(\gamma_{\min}(Z)/4\) is a universally valid inner
bound on the cell around the observed configuration. It is sufficient,
not generally maximal.

The factor four is not a removable artifact of loose notation. Under a
uniform per-label displacement budget, the focal label can move so as to
increase the selected distance and decrease the unselected distance,
while the two candidates move in the same adverse directions. These four
contributions can exhaust the original gap. Tighter certificates must
therefore exploit additional structure, such as label-specific motion
bounds, constrained perturbation directions, correlations among
displacements, or direct configuration-specific optimization.

\hypertarget{meaning-of-the-taylorgreen-application}{%
\subsection{7.2 Meaning of the Taylor--Green
application}\label{meaning-of-the-taylorgreen-application}}

For the frozen Taylor--Green information coordinates, the certified
radius proves that sufficiently small simultaneous coordinate
perturbations cannot alter any directed 5-NN membership at the three
included time strata. It consequently preserves the symmetrized graphs,
local clustering coefficients, frozen percentile rule, strict-core
labels, exposures, and observational next-time entry indicators defined
from those discrete inputs.

The numerical magnitude of the radius should be interpreted in the
coordinate system in which it was proved. The coordinates were globally
standardized and the normalization was held fixed. The certificate is
not automatically a radius in raw physical units, and it does not cover
a pipeline that re-estimates its scaling after perturbation. Its
principal value is logical and reproducible: it replaces an informal
assertion of robustness with a strict, second-implementation-checked
sufficient condition tied to frozen artifacts.

The certificate also does not turn the preserved structure into a
validated predictor. Structural invariance answers whether the same
analysis object is recovered under an allowed perturbation. Predictive
replication asks whether an association or phase pattern generalizes to
held-out data. The F32B-T temporal holdout did not support the complete
preregistered phase-switch pattern: boundary proximity retained a
favorable directional signal, but the full criterion failed, including a
coefficient with the wrong sign and a bootstrap interval for the AUC
difference crossing zero. This negative prospective result is compatible
with exact local stability of the graph and must not be softened by it.

\hypertarget{static-and-temporal-guarantees}{%
\subsection{7.3 Static and temporal
guarantees}\label{static-and-temporal-guarantees}}

The dynamic theorem supplies a clean bridge from spatial margin to time
only when the information-space trajectory has a bound valid throughout
the interval. If \(D_\infty(Z(t_0),Z(t_0+s))\le Ls\), then
\(\gamma_i(Z(t_0))/(4L)\) is a certified lower bound on the first possible
rewiring time for focal label \(i\). The label-specific result can improve
this bound by charging only the focal, selected, and unselected labels
relevant to each competing pair.

The F34 audit shows why this hypothesis must be treated as substantive.
Dense numerical resolution of the underlying flow and deformation
gradient does not certify dense regularity of a derived coordinate that
is rebuilt by clustering and boundary selection. An endpoint
displacement can miss arbitrarily faster intervening motion, while a
discrete change of defining branch can defeat a derivative argument that
ignores the switching mechanism. Reporting
\texttt{STOP\_NO\_INTERVAL\_VALID\_BOUND} therefore protects the
distinction between a theorem and an unsupported numerical substitution.

This STOP is informative rather than inconclusive. It identifies exactly
what a future specialization must supply: a fixed smooth information map
with a derivative bound, a certified piecewise-smooth description with
switching control, a validated dense-time reconstruction of the entire
map, or a modulus of continuity valid across combinatorial changes.
Until then, only the observed-time static certificates are authorized.

\hypertarget{broader-implications}{%
\subsection{7.4 Broader implications}\label{broader-implications}}

Rank-margin certification is applicable beyond the present flow example
whenever an analysis begins with a finite metric configuration and uses
a deterministic nearest-neighbor rule. Neighborhood graphs already
underpin geometric data-analysis methods {[}2,3{]}, while robustness has
been studied for k-NN classification under several perturbation models
{[}4--6{]}. The present theorem is deliberately agnostic about how the
points were generated; its strength and limitation both arise from that
abstraction.

The framework also suggests a useful reporting discipline. A study
should state the metric, labels, value of \(k\), tie-breaking convention,
normalization rule, perturbation norm, observed boundary margin,
strictness of the radius, and exact downstream dependency. For temporal
claims it should additionally state the interval and the source of the
interval-valid motion bound. These details make clear which parts of a
pipeline are mathematically protected and which remain empirical.

\hypertarget{next-directions}{%
\subsection{7.5 Next directions}\label{next-directions}}

Several extensions are natural. Heterogeneous perturbation budgets could
yield sharper focal-label guarantees than one global radius.
Configuration-specific optimization could estimate or certify the exact
distance to a neighbor-changing boundary under a chosen perturbation
model. Stability under recomputed normalization would require
propagating perturbations through the scaling map. When the minimum
boundary margin is zero, exact invariance may be replaced by bounds on
the number or weight of changing edges. Probabilistic measurement models
could convert deterministic margins into rewiring-risk bounds. Finally,
topological or persistence summaries could be studied across the
adjacent combinatorial cells rather than only within a single cell.
Persistence theory provides established stability results for
topological summaries under perturbation {[}7{]}, but transferring those
ideas to the present labeled, rank-defined graph cascade would require a
separate construction and theorem.

For the present application, the most direct route to a continuous-time
certificate is to redesign or augment the information map so that its
temporal regularity can itself be certified. Such work should be
preregistered as a new analysis and should preserve the existing F30
certificate and F34 STOP record unchanged.

\hypertarget{limitations}{%
\section{8. Limitations}\label{limitations}}

The general theorems concern finite labeled configurations, a fixed
metric, a fixed value of \(k\), and a fixed tie-breaking rule. They do not
cover label creation or loss, adaptive selection of \(k\), metric learning
performed after perturbation, or stochastic changes in the analysis rule
unless those operations are explicitly included in a new bound.

The radius \(\gamma_{\min}(Z)/4\) is a sufficient worst-case guarantee.
Except for the general sharpness statement, it is not claimed to equal
the maximal safe radius of the frozen Taylor--Green configuration. The
strict inequality must be retained because equality can produce a
boundary tie.

The numerical application is limited to one frozen construction, its
specified standardized coordinates, 585 labels, directed \(k=5\), and the
three observed starting-time strata. Recomputed normalization,
raw-coordinate perturbations, other values of \(k\), other samples, and
other flows require separate certificates. The arithmetic verified by a
second implementation within the same research workflow establishes the
stated finite rank margin; it does not
establish physical error bounds for every upstream simulation or
measurement process.

The continuous-time specialization is not numerically available. The
current artifacts do not establish interval-wide regularity, positivity,
or a motion bound for the full map \([d_B,\log A_B]\). In particular,
dense integration of the physical trajectories and \(J_F\) does not
control intervening changes in cluster assignments, boundary edges,
boundary normals, or \(d_B\). Secant slopes between sampled times cannot
replace the missing interval-valid bound.

The deterministic cascade result protects only objects that are fixed
functions of the unchanged labeled neighbor family and frozen
conventions. It does not protect quantities that use additional
perturbed inputs, retrained parameters, stochastic algorithms, or
altered thresholds. Nor does it establish that a preserved observational
outcome corresponds to an invariant physical event.

No claim is made here about causality, universal precedence, predictive
accuracy, stability of Navier--Stokes solutions, or the same radius
across dynamical regimes. The complete prospective phase-switch pattern
was not reproduced in F32B-T and is not presented as supported. The
reference set in this working draft establishes the immediate
methodological context but cannot prove universal priority; the
manuscript therefore uses ``we establish'' rather than ``for the first
time.''

\hypertarget{conclusion}{%
\section{9. Conclusion}\label{conclusion}}

A positive directed \(k\)-nearest-neighbor boundary margin gives an
explicit and sharp worst-case local guarantee for finite labeled
configurations: simultaneous per-label perturbations smaller than one
quarter of the margin preserve the directed neighbor family, and hence
every fixed deterministic construction based solely on it. With an
interval-valid information-space motion bound, the same result yields a
certified lower bound on the first possible rewiring time.

For the frozen Taylor--Green future-information graphs, exact interval
arithmetic and a separately implemented checker within the same research
workflow certified a positive static perturbation radius at three observed
time strata. This establishes
local combinatorial invariance in the frozen standardized coordinates.
The available artifacts did not establish the additional temporal
regularity needed for a numerical continuous-time horizon, so the
application remains correctly classified as \texttt{DISCRETE\_ONLY} with
\texttt{STOP\_NO\_INTERVAL\_VALID\_BOUND} for that extension.

The resulting conclusion is intentionally bounded: finite structural
stability has been proved and certified for the stated configuration,
whereas temporal specialization and predictive generalization require
separate evidence. Keeping those layers distinct turns a potentially
broad robustness claim into a precise, reproducible, and falsifiable
one.

\hypertarget{declarations}{%
\section{Declarations}\label{declarations}}

\textbf{Competing interests.} The author declares no competing
interests.

\textbf{Funding.} No external funding was received for this study.

\textbf{Ethics approval.} Not applicable. The study reports mathematical
results and a computational fluid-dynamics application; it includes no
human participants, identifiable human data, or animal experiments.

\textbf{Data and code availability.} The frozen numerical artifacts,
second-implementation checker outputs, and source files underlying the
Taylor--Green certificate are maintained with cryptographic checksums. They are not yet publicly archived. A public archival repository with a
persistent identifier is planned for a subsequent version of this manuscript.

\textbf{Author contributions.} H.S. conceived the study, developed the
future-information application, interpreted the results, and takes
responsibility for the manuscript. Computational and language tools,
including AI-assisted drafting and code review, were used under the author's
supervision. The author reviewed the mathematical statements, numerical
results, and final manuscript and assumes full responsibility for the work.

\hypertarget{references}{%
\section{References}\label{references}}

\begin{enumerate}
\def\labelenumi{\arabic{enumi}.}
\tightlist
\item
  T. M. Cover and P. E. Hart, ``Nearest Neighbor Pattern
  Classification,'' \emph{IEEE Transactions on Information Theory},
  13(1), 21--27, 1967. https://doi.org/10.1109/TIT.1967.1053964
\item
  J. B. Tenenbaum, V. de Silva, and J. C. Langford, ``A Global Geometric
  Framework for Nonlinear Dimensionality Reduction,'' \emph{Science},
  290(5500), 2319--2323, 2000.
  https://doi.org/10.1126/science.290.5500.2319
\item
  M. Belkin and P. Niyogi, ``Laplacian Eigenmaps for Dimensionality
  Reduction and Data Representation,'' \emph{Neural Computation}, 15(6),
  1373--1396, 2003. https://doi.org/10.1162/089976603321780317
\item
  Y. Wang, S. Jha, and K. Chaudhuri, ``Analyzing the Robustness of
  Nearest Neighbors to Adversarial Examples,'' in \emph{Proceedings of
  the 35th International Conference on Machine Learning}, Proceedings of
  Machine Learning Research, vol.~80, pp.~5133--5142, 2018.
  https://proceedings.mlr.press/v80/wang18c.html
\item
  A. Z. Fan and P. Koutris, ``Certifiable Robustness for Nearest
  Neighbor Classifiers,'' arXiv:2201.04770, 2022.
  https://doi.org/10.48550/arXiv.2201.04770
\item
  N. Fassina, F. Ranzato, and M. Zanella, ``Robustness Verification of
  k-Nearest Neighbors by Abstract Interpretation,'' \emph{Knowledge and
  Information Systems}, 66(8), 4825--4859, 2024.
  https://doi.org/10.1007/s10115-024-02108-4
\item
  D. Cohen-Steiner, H. Edelsbrunner, and J. Harer, ``Stability of
  Persistence Diagrams,'' \emph{Discrete \& Computational Geometry}, 37,
  103--120, 2007. https://doi.org/10.1007/s00454-006-1276-5
\end{enumerate}

\end{document}